\documentclass[11pt]{amsart}
\usepackage{amsmath,amssymb,enumerate}
\usepackage{epsf,epsfig,amsfonts,graphicx,color}  
 \usepackage{pgf,tikz}
 \usetikzlibrary{arrows}

 \date{September 11, 2014}   
  
 \numberwithin{equation}{section}   

\renewcommand{\r}{\mathbb{R}}

\newtheorem{theorem}{\rm\bf Theorem}[section]
\newtheorem{proposition}[theorem]{\rm\bf Proposition}
\newtheorem{lemma}[theorem]{\rm\bf Lemma}
\newtheorem{corollary}[theorem]{\rm\bf Corollary}
\theoremstyle{definition}

\newtheorem{definitions}[theorem]{\rm\bf Definitions}

\newcommand{\weg}[1]{}

\title{Isometries of two dimensional Hilbert Geometries}

\author{Vladimir S. Matveev} 
\address{Institut f\"ur Mathematik, Friedrich-Schiller Universit\"at Jena\\
07737 Jena, Germany}  
\email{vladimir.matveev@uni-jena.de}

\author{Marc Troyanov} 
\address{Section de Math{\'e}matiques,  
\'Ecole Polytechnique F{\'e}derale de Lausanne, station 8,
1015 Lausanne - Switzerland} 
\email{marc.troyanov@epfl.ch}

\begin{document}
\bigskip

\begin{abstract} 
We prove that any  isometry  between  two dimensional Hilbert geometries   is a projective transformation unless the domains are interiors of triangles.
\medskip
 
\end{abstract}

\thanks{We thank the Friedrich-Schiller-Universit\"at Jena, EPFL and the Swiss national Science Foundation for their support.}
 
\maketitle

 \vspace{-0.3cm}

\hfill {\small \it  Dedicated to Pierre de la Harpe on his seventieth birthday.}

 \vspace{0.3cm}

\section{Introduction}

\begin{figure}
\centering
\begin{tikzpicture}[line cap=round,line join=round,>=triangle 45,x=0.6cm,y=0.5cm]
\clip(-0.3,-3) rectangle (8.8,4.5);
 \draw  [fill= gray!20] plot  [smooth cycle]  coordinates  { (0.24,-0.54)  
 (1.3,-1.79)  
 (4.64,-2.69) 
 (7,-2.38)  
 (8.18,-0.46)  
 (5.77,2.36)  
 (3.56,3.1 )  
 (1.18,2.61)   
 (0.19,0.93)  }; 
\draw (0.19,0.93)-- (8.18,-0.46);
\fill   (8.18,-0.46) circle (1.5pt);
\draw  (8.47,-0.24) node {$\bar y$};
\fill   (0.19,0.93) circle (1.5pt);
\draw  (-0.1,1.08) node {$\bar x$}; 
\fill   (2.33,0.56) circle (1.5pt);
\draw  (2.33,0.93) node {$x$};
\fill   (5.52,0) circle (1.5pt);
\draw  (5.67,0.42) node {$y$};
\end{tikzpicture}
 \caption{The points $\bar x$ and $\bar y$}\label{def_hilbert}
\end{figure}
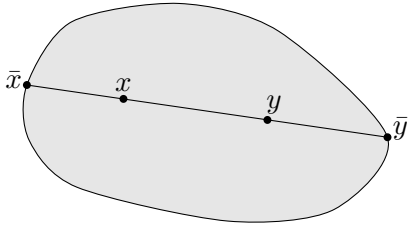

The  \emph{Hilbert distance} between two points $x$ and $y$ in a bounded convex domain $\Omega$ of  $\mathbb{R}^n$  is
defined as
\begin{equation} \label{formula_hilbert} 
  d(x,y):= \ln\bigl( (x,y;\bar x,\bar y)\bigr) := \ln\left( \frac{|\bar y - x|}{|\bar y -y|} : \frac{|\bar x - x|}{|\bar x -y|}  \right), 
\end{equation}
where $|u-v|$ denotes the usual Euclidean length between two points $u$ and $v$ in  $\mathbb{R}^n$,  
and $\bar x$ and $\bar y$ are as on Fig. \ref{def_hilbert}.   It is well known,  that the distance function  $d$
satisfies the standard requirements of a distance function, the only nontrivial  point to check being the triangle inequality,
see for example \cite{Hilbert} or  \cite[\S 1]{Harpe}. This distance has been introduced by Hilbert in \cite{Hilbert} and we refer to \cite{Handbook} for a presentation
of both classic and contemporary aspects of  Hilbert geometry\footnote{Although we have assumed that  $\Omega \subset \r^n$ is a bounded convex domain, the Hilbert distance (\ref{formula_hilbert}) is well defined for the
more general class of \emph{proper} convex domains.  A convex domain  in  $\mathbb{RP}^n$ is proper if it does not contain any full affine line. 
It is known that a convex domain is proper if and only if it is projectively equivalent to a bounded convex domain.  For convenience we will therefore consider only bounded convex domains.}.  
   
Recall that straight lines, convexity, and the cross ratio of four aligned points are  invariant under projective transformations, this implies immediately
that if $f : \mathbb{RP}^n \to  \mathbb{RP}^n$ is a projective transformation, then its restriction to $\Omega$ defines an isometry $f : \Omega \to f(\Omega)$,
with respect to the Hilbert distances in $\Omega$ and $f(\Omega)$.  (We consider $\mathbb{R}^n$ as a subset of $\mathbb{RP}^n$  by identifying it with an  affine chart, the Hilbert metric inside $\Omega$  does not depend on the choice of the affine chart.)
The converse to this statement is not always true: some special Hilbert geometries admit isometries which are not projective   transformations.
The simplest example is given by the simplex and is discussed in details in dimension 2 by Pierre de la Harpe in \cite{Harpe}. This author asked
for a    full description of all isometries in  Hilbert geometry and a complete answer in finite dimension  has recently been
obtained by Cormac Walsh in \cite{Walsh}. Note also that the same author, together with Bas Lemmens, previously described all isometries of
polyhedral Hilbert geometries  in \cite{LemmensWalsh2011}, while Bas Lemmens, Mark Roelands and  Marten Wortel gave some partial 
results in infinite dimension  in \cite{LRW}.

Our goal in this paper is to give a  short proof of the following two dimensional result:

\begin{theorem}\label{main} 
Let $\Omega_1$ and $\Omega_2$ be two bounded convex domains in the plane $\mathbb{R}^2$ and 
$d_1, d_2$ be the corresponding  Hilbert metrics.  Suppose that $\Omega_1$ is not  the interior of
a triangle, then every isometry   $f:(\Omega_1, d_1) \to (\Omega_2, d_2)$ 
is the restriction of a   projective transformation of $\mathbb{RP}^2$. 
\end{theorem}

As mentioned above, this result is false if $\Omega_1$ is  the interior of a triangle. In that case $(\Omega_1, d_1)$ is 
isometric to a Minkowski plane whose unit ball is a regular  hexagon and its group of isometries is not difficult to describe,
see \cite{Harpe}.  Recall also that above theorem is a special case of the result of C. Walsh \cite[Theorem 1.3]{Walsh}. For the
case of quadrilaterals, the result is also proved by P. de la Harpe in \cite[Proposition 4]{Harpe}.

Our proof uses completely different methods from those in Walsh's paper. It is quite direct and  only based on the description of 
metric geodesics in Hilbert geometry, together with a quite old and nontrivial result from line geometry which is due W. Prenowitz.

\section{The case of strictly convex domains}

It will be convenient to start with the case of a strictly convex domain. In fact we will prove the following 
result:

\begin{proposition}\label{prop.strict}
Assume that $\Omega_1\subset \r^n$ is a strictly convex domain,   then every isometry   $f:(\Omega_1, d_1) \to (\Omega_2, d_2)$ 
is the restriction of a   projective transformation of $\mathbb{RP}^n$. 
\end{proposition}

This result is proved in \cite[Proposition 3]{Harpe}, but we shall give a slightly more direct proof. The result has recently been 
extended in infinite dimension in \cite[Theorem 1.2]{LRW}.

The proof is based on the structure of geodesics for the Hilbert distance. It is easy to check from the definition of the Hilbert distance that if three points $x,y,z \in \Omega$ are aligned and  $z \in [x,y]$, then $d_1(x,y) = d_1(x,z)+d_1(z,y)$. In other words the intersection of Euclidean straight lines with $\Omega_1$ are
geodesics for the Hilbert metric.  Furthermore, the following fact is classical  (see  \cite[Proposition 2]{Harpe} or \cite[Theorem 12.5]{PT_Funk}):

\begin{lemma}\label{lem.uniquegeodesic}
Let $p$ and $q$ be two points on the boundary of $\Omega_1$, and suppose that at least one of them is an extreme point of $\Omega_1$.
Then the open interval $(p,q)$ is the unique geodesic between any pair of its point, that is if $x,y \in (p,q)$ and $z \in \Omega_1$, then 
$d_1(x,y) = d_1(x,z)+d_1(z,y)$ if and only if $z \in [x,y]$.
\end{lemma}

\medskip

\textbf{Proof of Proposition \ref{prop.strict}.}   Let $f: \Omega_1 \to \Omega_2 $ be an isometry for the Hilbert distances, where 
 $\Omega_1\subset \r^n$ is strictly convex. From the previous Lemma,  it then follows that  the affine segment $[x,y]$ between two points
 $x,y \in \Omega_1$ is the unique geodesic joining these two points. Since $f$ is an isometry, there is also a unique geodesic joining
 the images $f(x)$ and $f(y)$ in $\Omega_2$ and because the Euclidean segment $[f(x),f(y)] \subset \Omega_2$ is known to be geodesic
 we conclude that $f$ maps the segment $[x,y]\subset \Omega_1$ to the segment $[f(x),f(y)] \subset \Omega_2$. Since $x$ and $y$ are
 arbitrary points in $\Omega_1$, we conclude that $f$ is a \emph{local collineation}, that is a mapping sending  Euclidean segments to  
 Euclidean segments.
 The conclusion now follows from the local version of the fundamental theorem of projective geometry (see e.g. \cite[Lemma 4]{Shiffman}), 
 which states that any local collineation defined in some open set of the  real projective space $\mathbb{RP}^n$ is the restriction of a projective transformation.
 
 \qed

\section{Proof of the main  Theorem}

The proof of  Theorem  \ref{main} will be based on a 1935  result of Prenowitz   \cite{Prenowitz} which generalizes  the fundamental theorem of projective geometry 
in dimension 2.  We will need the following definitions.
   
\begin{definitions}
Let $U$ be a plane domain, that is an open connected nonempty subset of  $\mathbb{R}^2$. By a \emph{line in $U$} we mean a connected component of  the intersection 
of a Euclidean straight line with $U$.  A   \emph{family of lines}  in $U$  is a partition of $U$ by lines, that is   a collection of  lines in $U$  such that 
 each point of $U$ lies on exactly  one line of the collection.  If all lines in a family  concourse  to a common point $A$, the family is called a \emph{pencil  with pole $A$}.   A (linear) \emph{$n$-web} in $U$ is a set of $n$  families of lines on $U$ such that no two families  have a common line.
 \end{definitions}
  
 \smallskip
 
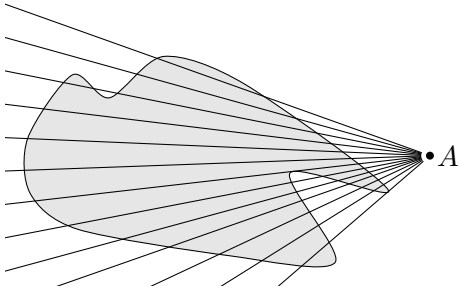
\begin{figure}
\centering
\begin{tikzpicture}[line cap=round,line join=round,>=triangle 45,x=0.6cm,y=0.5cm]
\clip(-0.3,-3) rectangle (10,4.5);
 \draw  [fill= gray!20] plot  [smooth cycle]  coordinates  { (0.24,-0.54)  
 (1.3,-1.5)  
 (4.64,-2.2) 
 (7,-2.38)  (6,0)
 (8.18,-0.46)  
 (5,2.36)  
 (3.2,3.1)  (2,2)
 (1.18,2.61)   
 (0.19,0.93)  }; 
\draw (9.1,0.46) -- (-1.5,-3);
\draw (9.1,0.46) -- (-1.5,-2);
\draw (9.1,0.46) -- (-1.5,-1);
\draw (9.1,0.46) -- (-1.5,0);
\draw (9.1,0.46) -- (-1.5,3);
\draw (9.1,0.46) -- (-1.5,2);
\draw (9.1,0.46) -- (-1.5,1);
\draw (9.1,0.46) -- (-1.5,4);
\draw (9.1,0.46) -- (-1.5,-4);
\draw (9.1,0.46) -- (-1.5,-5.5);
\draw (9.1,0.46) -- (-1.5,5.0);
\draw (9.1,0.46) -- (-1.5,-7.5);
\draw (9.1,0.46) -- (-1.5,-10.5);
\fill  [color=white] (9.1,0.46) circle (3.3pt);
\fill  [color=black] (9.1,0.46) circle (1.5pt);
\draw  (9.5,0.45) node {$A$};
\end{tikzpicture}
\caption{A  pencil of lines covering a plane domain.}\label{family1}
\end{figure}

Figure \ref{family1} shows a pencil with pole $A$ in the domain $U$. 
By taking  the pencils through   $n$ pairwise distinct poles $A_1, \dots A_n \not\in U$ we obtain
an $n$-web in any subdomain $U' \subset U$ disjoint from any line through a pair of distinct points $A_i$,$A_j$.

\begin{theorem}[Prenowitz 1935] \label{thm.Prenowitz} 
A one to one continuous map defined in a plane domain 
that  carries a 4-web  into a 4-web is the restriction of a  projective transformation.
\end{theorem}  
 
 Recall that by Brouwer's theorem, an injective continuous map defined in a  domain of $\r^n$ is a homeomorphism onto
 its image. The above result is proved in \cite{Prenowitz}; a much simpler proof is given in \cite{Kasner} assuming the map
 is a diffeomorphism.  Some generalization in higher dimensions are given in \cite{AAS}.   
 
 The following corollary will  be useful in the proof:

\begin{corollary} \label{cor.Prenowitz} 
Let $f : U \to \mathbb{R}^2$  be a one to one continuous map defined in a domain 
$U\subset \mathbb{R}^2$ and let $A_1,\cdots ,A_5 \in \mathbb{R}^2$
be five pairwise distinct points. Assume that $f$ maps the intersection of any line through 
$A_j$ with $U$ to a straight line ($1 \leq j \leq 5$).
Then $f$ is the restriction of a projective transformation.
\end{corollary}

\begin{figure}[ht] 
\begin{tikzpicture}[line cap=round,line join=round,>=triangle 45,x=1.3cm,y=1.3cm]
\clip(-.5,-.5) rectangle (7,7);
\fill[color=gray!20,fill=gray!20]  (2,4.5) -- (1.39,4.02) -- (1.5,3) -- (2.6,2.28) -- (3.83,2.7) -- (4.77,3.5) -- (4.59,4.76) -- (3.64,4.37) -- cycle;
\draw [color=gray] (0,1.5)-- (3.84,5.48);
\draw [color=gray] (0,1.5)-- (4.62,5.26);
\draw [color=gray] (0,1.5)-- (5.22,4.55);
\draw [color=gray] (0,1.5)-- (5.27,3.29);
\draw [color=gray] (0.11,5.75)-- (3.21,1.18);
\draw [color=gray] (0.11,5.75)-- (3.45,1.85);
\draw [color=gray] (0.11,5.75)-- (4.07,2.17);
\draw [color=gray] (0.11,5.75)-- (4.88,2.49);
\draw [color=gray] (0.11,5.75)-- (5.46,2.63);
\draw [color=gray] (0.11,5.75)-- (5.34,3.55);
\draw [color=gray] (0.11,5.75)-- (5.46,4.18);
\draw [color=gray] (2.5,0.5)-- (1.13,4.95);
\draw [color=gray] (2.5,0.5)-- (1.69,4.91);
\draw [color=gray] (2.5,0.5)-- (2.08,5.06);
\draw [color=gray] (2.5,0.5)-- (2.53,5.34);
\draw [color=gray] (2.5,0.5)-- (3.08,5.25);
\draw [color=gray] (2.5,0.5)-- (3.82,5.04);
\draw [color=gray] (2.5,0.5)-- (4.94,4.77);
\draw [color=gray] (2.5,0.5)-- (4.45,4.91);
\draw [color=gray] (6,3)-- (1.5,2.5);
\draw [color=gray] (6,3)-- (1.08,3.42);
\draw [color=gray] (6,3)-- (1.5,4.3);
\draw [color=gray] (6,3)-- (2.08,5.06);
\draw [color=gray] (6,3)-- (3.18,5.04);
\draw [color=gray] (4,6)-- (1.2,3.64);
\draw [color=gray] (4,6)-- (1.5,2.5);
\draw [color=gray] (4,6)-- (2.02,2.34);
\draw [color=gray] (3.84,5.48)-- (2.89,2.23);
\draw [color=gray] (4,6)-- (3.43,2.33);
\draw [color=gray] (4,6)-- (4.21,2.63);
\fill  (0,1.5) circle (1.5pt);
\draw (4.57,3.6)  node {$U$} ;
\draw (-0.2,1.52) node {$A_1$};
\fill (2.5,0.5) circle (1.5pt);
\draw (2.6,0.32) node {$A_2$};
\fill (6,3) circle (1.5pt);
\draw (6.29,2.96) node {$A_3$};
\fill (4,6) circle (1.5pt);
\draw (4.15,6.23) node {$A_4$};
\fill (0.11,5.75) circle (1.5pt);
\draw (-0.1,5.9) node {$A_5$};
\end{tikzpicture}
\caption{A polygonal region $U$ covered by 5 pencils. Corollary \ref{cor.Prenowitz} states that a homeomorphism defined in $U$ carrying all those lines
to some families of lines is a projective transformation.}\label{family}
\end{figure}
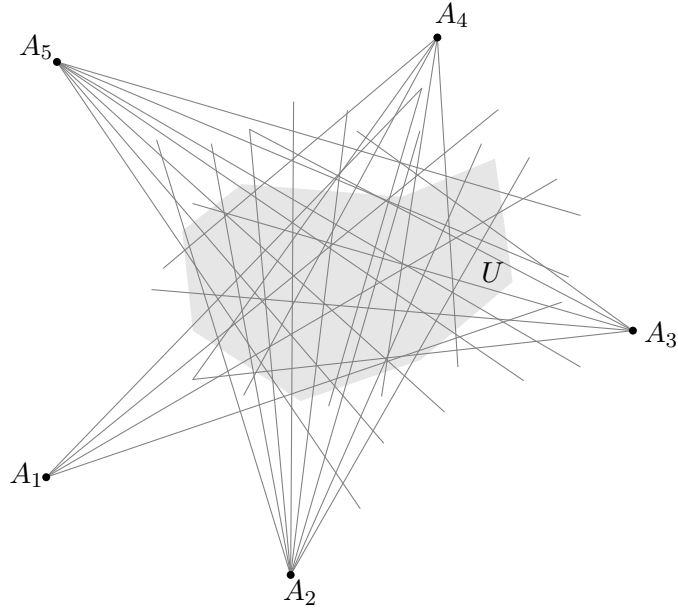

\begin{proof}
There are 10 lines through any pair of the points $A_j$   and the pairwise intersections of those 10 lines determine (at most)  20 points\footnote{10 distinct lines in a projective plane   define $\binom{10}{2}=45$ intersection points counted with multiplicity,  the 5 points $A_j$ have multiplicity $6$.}. 
Let us denote by  $\mathcal{I}$ this set and call it the set of \emph{intersection points}. 
For any point $X \in U\setminus \mathcal{I}$, at least four of the directions $\overrightarrow{XA_j}$ are mutually distinct  and this property holds
in a neighborhood $V$ of $X$. The pencils with poles the corresponding four points $A_j$ form a 4-web in  $V$, see Figure \ref{family},  which is mapped by $f$ to 
a 4-web in $f(V)$. By Theorem \ref{thm.Prenowitz},  we know that the restriction of $f$ to $V$ is the restriction of a  projective transformation.   
By real analyticity,  two projective transformations that coincide on an open subset coincide everywhere. Since $U\setminus \mathcal{I}$
is connected the restriction of  $f$ to $U\setminus \mathcal{I}$ is a projective transformation and since $\mathcal{I}$ is finite, 
$f$ is a projective transformation on the whole domaine $U$ by continuity.
\end{proof} 
 
  \textbf{Proof of  Theorem \ref{main}. }
Recall that we assumed that the bounded convex domain $\Omega_1 \in \mathbb{R}^2$ is not the interior of a triangle.
We first assume that $\Omega_1$ is also not a quadrilateral. 
Then, the boundary   $\partial\Omega_1$ contains at least five distinct  extremal points $A_1,A_2,A_3,A_4,A_5$. 
Because the points $A_j$ are extreme points of $\Omega_1$,  Lemma \ref{lem.uniquegeodesic} implies that  each line 
through one of the point $A_j$  intersects $\Omega_1$ on a unique geodesic (for the Hilbert distance) between any of its pair of point. 
Since $f$ is an isometry, it sends each line from the five pencils  to a straight line in $\Omega_2$
and it follows from Corollary \ref{cor.Prenowitz} that $f$ is the restriction of a projective transformation.

\smallskip

 Suppose now that $\Omega_1$ is a quadrilateral with vertices $ABCD$. The vertices are extreme points of $\Omega_1$,
 therefore, by Lemma \ref{lem.uniquegeodesic}, any line through a vertex defines a unique geodesic for the Hilbert distance and it is thus mapped on a line
 by the isometry $f$. The pencils with poles the four vertices form a 4-web in each connected component of the complement of the diagonals (these components
 are the triangles  and $DAM$, where $M$ is the intersection of the diagonals).   
 Using  Prenowitz' Theorem  \ref{thm.Prenowitz}, we conclude that   the restriction of $f$ to each of the triangles $ABM, \  BCM, \ CDM, \ DAM$   is a projective transformation.

 \begin{figure}
\centering
\definecolor{uuu}{rgb}{0.9,09,09}
\begin{tikzpicture}[line cap=round,line join=round,>=triangle 45,x=0.7cm,y=0.7cm]
\clip(-0.5,-0.5) rectangle (4.5,4.5);
\fill[color=gray!20] (0,0) -- (4,0) -- (4,4) -- (0,4) -- cycle;
\fill[color=gray!45] (0.2,0.07) -- (2,1.9) -- (3.8,0.07) -- cycle;
\fill[color=gray!45]  (3.93,0.2) --(2.1,2) -- (3.93,3.8)-- cycle;
\draw (0,4)-- (4,0); 
\draw (0,0)-- (4,4); 
\begin{scriptsize}
\fill  (0,0) circle (1.5pt);
\draw (-0.23,-0.02) node {$A$};
\fill [color=black] (4,0) circle (1.5pt);
\draw (4.21,-0.06) node {$B$};
\fill [color=black] (4,4) circle (1.5pt);
\draw (4.23,4.17) node {$C$};
\fill [color=black] (0,4) circle (1.5pt);
\draw (-0.2,4.19) node {$D$};
\fill (2,2) circle (1.5pt);
\draw (1.55,2.08) node {$M$};
\end{scriptsize}
\end{tikzpicture}
 \caption{\small The restriction of $f$ to the ``dark-gray'' triangles  $ABM$ and $BCM$ is a projective transformation, and $f$ sends $AC$ to a straight line. Then, the image of $ABC$ is a triangle and  the restriction of $f$ to it 
 is a projective transformation.}  
 \label{square}
\end{figure}
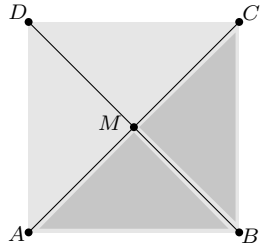

Consider two adjacent  such triangles,   and consider the  $f$-image of their union, see Fig \ref{square}.  Since the restriction of $f$ to each of these triangles is a projective transformation, the image of its union is two triangles. By continuity it has a common edge. Since the image of the line $AC$ is a straight line, the closure of the image of the union of these triangles  is a triangle. Furthermore, the map $f$ sends any line through
$A$ or $B$ to a line, we thus conclude that $f$  restricted to  the triangle $ABC$ is a projective transformation (see also the Corollary
in \cite{Prenowitz}  page 567). Similarly, the restrictions  of $f$ to $BCD$ and to $CDA$   are projective transformations, which implies that the map $f$ on the whole quadrilateral $ABCD$ is the restriction of  projective transformation as  desired. 

 \qed

\medskip

The first author would like to thank Richard   Schwartz for attracting the interest to this problem and  both authors thank  Cormac Walsh,  Bas Lemmens, Patrick Foulon and  Inkang  Kim 
for useful discussions on the subject of this paper.



\begin{thebibliography}{999}

\bibitem{AAS}  S.  Artstein-Avidan and B. Slomka, 
\emph{The Fundamental Theorems of Affine and Projective Geometry Revisited.}
preprint (work in preparation).
 
\bibitem{Busemann1944} H. Busemann, \emph{Local metric geometry. } {Trans. Amer. Math. Soc.}  {\bf 56}(1944) 200--274.
   
\bibitem{Busemann1953} H. Busemann and P. J. Kelly, \emph{Projective geometry and projective metrics},  Academic Press  (1953).

\bibitem{Busemann1955} H. Busemann,  \emph{The geometry of geodesics}, Academic Press  (1955), reprinted by Dover in 2005.
 
\bibitem{Harpe} P. de la Harpe, \emph{
On Hilbert's metric for simplices.}  Geometric group theory, Vol. {\bf 1} (Sussex, 1991), 97--119,
London Math. Soc. Lecture Note Ser.,  {\bf 181}, Cambridge Univ. Press, Cambridge, 1993. 
 
\bibitem{Handbook}  \emph{Handbook of Hilbert geometry}, (ed. A. Papadopoulos and M. Troyanov), European Mathematical Society, Z\"urich, 2014.

\bibitem{Hilbert} D. Hilbert,  \emph{Ueber die gerade Linie als k\"urzeste Verbindung zweier Punkte.} Math. Ann. {\bf 46}(1895) 91--96.

 \bibitem{Kasner} E. Kasner,  
\emph{The characterization of collineations.} 
Bull. Amer. Math. Soc. {\bf 9}(1903), no. 10, 545--546. 

\bibitem{LRW} 
B. Lemmens, M. Roelands, M. Wortel, 
\emph{Isometries of infinite dimensional Hilbert geometries.}
arXiv:1405.4147

\bibitem{LemmensWalsh2011} B. Lemmens and C. Walsh,  
\emph{Isometries of polyhedral Hilbert geometries.}  
J. Topol. Anal. {\bf 3}(2011), no. 2, 213--241. 


\bibitem{PT_Funk} A. Papadopoulos and M. Troyanov, \emph{From Funk to Hilbert Geometry.} In: \emph{Handbook of Hilbert geometry}, 
(ed.   A. Papadopoulos and M. Troyanov), European Mathematical Society, Z\"urich, 2014.

 
\bibitem{Prenowitz}  W. Prenowitz,  
\emph{The characterization of plane collineations in terms of homologous families of lines.} Trans. Amer. Math. Soc. {\bf 38}(1935)  564--599.

\bibitem{Shiffman}  B. Shiffman, 
\emph{Synthetic projective geometry and Poincar\'e theorem on automorphisms of the ball. }
Enseign. Math. (2) {\bf 41}(1995), no. 3-4, 201--215. 

\bibitem{Walsh} C.  Walsh, \emph{ Gauge-reversing maps on cones, and Hilbert and Thompson isometries.}  	arXiv:1312.7871
   
\end{thebibliography}
\end{document}